\documentclass[letterpaper, conference]{IEEEtran} 

\IEEEoverridecommandlockouts 

\usepackage{amsmath}
\usepackage{amsfonts}
\usepackage{mathtools}
\usepackage{array}
\usepackage{gensymb}
\usepackage{graphicx}
\usepackage[export]{adjustbox}
\usepackage{caption}
\usepackage{subcaption}
\usepackage{verbatim}
\usepackage{multirow}
\usepackage{tikz}
\usepackage{eurosym}
\usepackage{enumerate}
\usepackage{amsthm}
\usetikzlibrary{shapes,arrows}
\usepackage{algorithm}
\usepackage{url}
\usepackage{algpseudocode}
\usepackage{accents}
\usepackage[normalem]{ulem}
\newcommand{\ubar}[1]{\underaccent{\bar}{#1}}

{{\slshape}\newtheorem{theorem}{Theorem}[section]}
{{\slshape}}
{{\slshape}}
{{\slshape}}
{{\slshape}}
{{\upshape}}
{{\upshape}}
{{\upshape}\newtheorem{remark}[theorem]{Remark}}
{{\upshape}}

\newcommand{\dual}{\text{dual}}
\newcommand{\prim}{\text{prim}}

\newcommand{\mat}[1]{\begin{bmatrix}#1\end{bmatrix}}
\newcommand{\costFcn}[3]{\min_{#1} ~~ & #2 \\ \text{s.t.} ~~ & #3}

\newcommand{\solver}{superADMM} 
\newcommand{\Solver}{SuperADMM}
\newcommand{\kpp}{{k+1}}

\title{\huge\bf \Solver: Solving Quadratic Programs Faster with Dynamic Weighting ADMM\thanks{This research was performed within the framework of the research program AquaConnect, funded by the Dutch Research Council (NWO, grant-ID P19-45) and public and private partners of the AquaConnect consortium and coordinated by Wageningen University and Research.}}

\author{P. C. N. Verheijen$^{1}$, D. Goswami$^{1}$ and M. Lazar$^{1}$%
\thanks{$^{1}$Department of Electrical Engineering, Eindhoven University of Technology, The Netherlands: 
\texttt{p.c.n.verheijen@tue.nl, d.goswami@tue.nl, m.lazar@tue.nl}}
} 

\begin{document}
\maketitle
\thispagestyle{empty}
\pagestyle{empty}

\begin{abstract}%
In this paper we develop an accelerated Alternating Direction Method of Multipliers (ADMM) algorithm for solving quadratic programs called \solver.
Unlike standard ADMM QP solvers, \solver\ uses a novel dynamic weighting method that penalizes each constraint individually and performs weight updates at every ADMM iteration. 
We provide a numerical stability analysis, methods for parameter selection and infeasibility detection.
The algorithm is implemented in \texttt{c} with efficient linear algebra packages to provide a short execution time and allows calling \solver\ from popular languages such as MATLAB and Python. A comparison of \solver\ with state-of-the-art ADMM solvers and widely used commercial solvers showcases the efficiency and accuracy of the developed solver.
\end{abstract}

\section{Introduction}
Quadratic Programming (QP) refers to a constrained optimization problem that minimizes a convex quadratic cost function subject to linear constraints. Quadratic programming is one of the most efficient approaches to solving various optimization problems, ranging from data fitting, static optimization, or constrained linear Model Predictive Control (MPC) \cite{ConvexOptBook, RawlingsBook}. Moreover, complex problems like nonlinear and mixed-integer programs can often be solved by successively approximating them with QPs, e.g., via Sequential Quadratic Programming (SQP) \cite{SQP} or Branch and Bound.

The development of quadratic programming solvers has been an active research topic for the last 70 years \cite{KKT:Kuhn1951}. The widely used approaches include active set solvers such as qpOASES \cite{qpOASES:Ferreau2014}, interior point methods such as Gurobi \cite{gurobi} or Mosek \cite{mosek}, gradient methods (cf. \cite{GradientMethod:Necoara2018}), and Augmented Lagrangian methods such as OSQP \cite{OSQP:Stellato2020} or QPALM \cite{hermans2019qpalm}. The continued search for faster QP solvers is necessary to facilitate real-time implementation of MPC \cite{InexactAugLag:Necoara2014}, as this ensures that the control input is computed within the sampling interval. More modern MPC applications, such as power electronics or autonomous driving, require a faster response time, while problems expand in size and complexity \cite{FastMPC:Wang2010}.

The aforementioned augmented Lagrangian methods are often solved using the Alternating Direction Method of Multipliers, or ADMM. ADMM decomposes problems into smaller subproblems that can often be solved without expensive Newton steps. As a result, the iterations are computationally inexpensive. As a drawback to being a first-order method, ADMM suffers from slow linear convergence near the optimum \cite{ConvADMM:Hong2017}\cite{LinConvDistGrad:Necoara2015}. However, if solution accuracy is not a high priority; the use of ADMM can outperform alternative QP solver methods. Unlike interior-point methods, ADMM can be warm-started, which is especially beneficial for the receding-horizon principle in MPC. However, its performance depends heavily on the choice of the penalty parameter in the augmented Lagrangian. To optimally select this penalty parameter, OSQP \cite{OSQP:Stellato2020} balances the primal and dual residuals. Alternatively, QPALM \cite{hermans2019qpalm} improves on this by assigning individual penalties per constraint and updating them conditionally, as follows. First, all constraints that are currently inactive are ignored in the problem (equivalently, setting their penalty to zero). Then, if the residual of the constraints does not decrease, the corresponding penalty increases. Regarding the parameter selection for ADMM methods, it has been proven for equality-constrained problems that if the penalty increases to infinity, the corresponding convergence rate will become \textit{superlinear} \cite{ADMMproof:Bertsekas1976} and as such, it is possible in principle to also achieve fast convergence close to the optimum. In practice, the increase of the penalty can remain bounded by a very large number instead of infinity.

In this paper, we introduce a new ADMM algorithm for solving QP problems, called \solver, with the following contributions:
\begin{itemize}
    \item a novel dynamic weighting method affecting the individual penalties on the constraints in the augmented Lagrangian;
    \item we provide a novel method that detects and retains numerical stability of the iterates of the solver;
    \item we provide an implementation in \texttt{c} for short execution time;
    \item evaluation on a set of benchmark problems, comparing \solver\ against the state-of-the-art in QP solving.
\end{itemize}

\Solver\ is a new ADMM-based solver for QPs with dynamic, per-constraint exponential penalty updates. With this novel dynamic weighting method, we penalize each constraint individually and exponentially increase all penalties if the corresponding constraint is on the active set, and exponentially decrease all others. Through this, the penalties approach infinity for all active constraints, which promotes superlinear convergence behavior for \solver\ for inequality constrained problems, hence the name. 
As a result, whereas standard ADMM often converges with a slow, linear convergence rate as it approaches the optimum, the convergence rate of \solver\ increases with each iteration. This results in significantly fewer iterations and yields a convergence from a low-accuracy solution to a high-accuracy solution in very few iterations. 
To maintain numerical stability, we introduce a bound on the penalty growth and decrease this bound if the numerical precision error exceeds the primal or dual residuals. This ensures practical convergence of the algorithm, and it is also used to terminate prematurely if a solution more accurate than the specified tolerances cannot be computed.
\Solver\ has been written in \texttt{c} and supports both dense and sparse problems. We compare the solver to state-of-the-art ADMM solvers such as OSQP and QPALM, and to the commercial solvers Gurobi and Mosek, for a set of different QP problem classes. The results show that \solver\ is competitive and can outperform all these solvers for most of the benchmarks considered in the literature \cite{OSQP:Stellato2020}. 

The remainder of this paper is structured as follows. In Section \ref{sec:prelim} we provide preliminaries on solving QPs using ADMM. We introduce the \solver\ algorithm in Section \ref{sec:superADMM}. In Section \ref{sec:implem} we provide some details about the implementation and some practical aspects of the solver. The benchmarks comparing the solver against the state-of-art are presented in Section \ref{sec:benchmarks} along with a compact description of OSQP and QPALM algorithms, which better highlights similarities and differences with \solver. Finally, we provide some concluding remarks in Section \ref{sec:conclusions}.



\section{ADMM preliminaries}\label{sec:prelim}
Consider a quadratic program of the following form:
\begin{equation}\label{eq:ADMM_prob}
    \begin{aligned}
        \costFcn{x}{\frac{1}{2}x^TPx + x^Tq} {l \leq Ax \leq u},
    \end{aligned}
\end{equation}
where $x\in\mathbb{R}^n$ is the vector of optimization variables, $P\in\mathbb{R}^{n\times n}$ is a positive semi-definite cost matrix and $q\in\mathbb{R}^n$ a linear cost vector. Both equality and inequality constraints are embedded in the upper-lower bound structure $l\leq Ax\leq u$, with $A\in\mathbb{R}^{m\times n}$ a linear mapping. This can be solved by replacing problem \eqref{eq:ADMM_prob} with an equality constrained problem, i.e.,
\begin{equation}\label{eq:ADMM_prob_z}
    \begin{aligned}
        \costFcn{x, z}{\frac{1}{2}x^TPx + x^Tq + g(z)} {Ax = z},
    \end{aligned}
\end{equation}
where the indicator function $g(z):\mathbb{R}^m \rightarrow \mathbb{R} \cup \{+\infty\}$ is defined as
\begin{equation}
    g(z) := \begin{cases}
        0, ~~ & \text{if}~ l\leq z \leq u \\
        \infty, ~~ &\text{else},
    \end{cases}
\end{equation}
and $z \in \mathbb{R}^m$ is an additional decision variable that directly operates on the constraint, i.e., $z=Ax$. If $l\leq z^* \leq u$, then it can be shown \cite{ADMMbook:Boyd2010} that the optimal solution of \eqref{eq:ADMM_prob_z} is the optimal solution of \eqref{eq:ADMM_prob}. The problem defined in \eqref{eq:ADMM_prob_z} can be solved using ADMM by constructing the augmented Lagrangian
\begin{equation}\label{eq:lagran}
    \begin{aligned}
        \mathcal{L}_\rho(x,z,y) :=& \frac{1}{2}x^TPx + x^Tq + g(z)+ y^T(Ax-z) \\ & + \frac{\rho}{2}\|Ax-z\|_2^2,
    \end{aligned}
\end{equation}
with dual variable $y\in\mathbb{R}^m$ and computing the following iterates until convergence
\begin{equation}\label{eq:ADMM_alg}
    \begin{aligned}
        x^{k+1} &:= \text{arg}\min_x \mathcal{L}_\rho(x, z^k, y^k) \\
        z^{k+1} &:= \text{arg}\min_z \mathcal{L}_\rho(x^{k+1}, z, y^k) \\
        y^{k+1} &:= y^k +\rho (Ax^{k+1}-z^{k+1}).
    \end{aligned}
\end{equation}
With respect to $g(z)$, the $z^{k+1}$ update in \eqref{eq:ADMM_alg} can be rewritten to 
\begin{equation}\label{eq:ADMM_algDirect}
    \begin{aligned}
        z^{k+1} &:= \Pi_{\left[l, u\right]}(Ax^{k+1} + \frac{1}{\rho}y^k),\\
    \end{aligned}
\end{equation}
where $\Pi_{\left[l,u\right]}(\cdot)$ is a projection operator, which clips the argument element-wise between $l$ and $u$. Convergence of the algorithm is tested by computing
\begin{equation}\label{eq:convergenceTerms}
\begin{aligned}
    r_\text{prim}^k &= \|Ax^k-z^k\|_\infty\\ 
	r_\text{dual}^k &= \|Px + q + A^Ty^k\|_\infty
   \end{aligned}
\end{equation}
and achieved if both
\begin{equation}\label{eq:convergenceCrit}
\begin{aligned}
    r_\text{prim}^k \leq \epsilon_\text{abs}, ~~ r_\text{dual}^k \leq \epsilon_\text{abs},
   \end{aligned}
\end{equation}
for some tolerance $\epsilon_\text{abs} > 0$. Finally, we extend the augmented Lagrangian with a proximal point term, i.e.,
\begin{equation}
    \begin{aligned}\label{eq:lagr_prox}
        \mathcal{L_{\rho,\sigma}}(x,x^k,z,y) := & \frac{1}{2}x^TPx + x^Tq + g(z)+ y^T(Ax-z) \\ & + \frac{\rho}{2}\|Ax-z\|_2^2 + \frac{\sigma}{2}\|x-x^k\|_2^2,
    \end{aligned}
\end{equation}
for some $\sigma\geq 0$. The convergence of the method with the proximal point is proven in \cite{ALM_Proximal:Rockafellar1976}. If $\sigma$ is kept appropriately small, the proximal point algorithm has little influence on the convergence rate. Furthermore, it has great numerical benefits, such as minimizing a strong convex function in the update of $x^\kpp$, as shown in \cite{ALM_Proximal:Rockafellar1976}. As a result, it is more robust against ill-conditioned problems.

\begin{remark}
    From this point, the \textbf{core} OSQP \cite{OSQP:Stellato2020} algorithm can be recovered by solving $\mathcal{L}_{\rho, \sigma}(x, x^k, z, y)$ in \eqref{eq:lagr_prox} using the (albeit slightly modified) ADMM steps in \eqref{eq:ADMM_alg}. To speed up convergence of the ADMM steps, OSQP updates the penalty parameter by balancing $r_\prim^k$ and $r_\dual^k$ if the gap between them is sufficiently large \cite[Sec.~5.2]{OSQP:Stellato2020}, i.e.,
    \begin{equation}\nonumber
        \begin{aligned}
            \rho^{k+1} = \rho^k \sqrt{\frac{r^k_\prim}{r_\dual^k}\frac{\mathrm{max}(\|Px^k\|_\infty, \|A^Ty^k\|_\infty, \|q\|_\infty)}{\mathrm{max}(\|Ax^k\|_\infty, \|z^k\|_\infty)}}.
        \end{aligned}
    \end{equation}
    The complete OSQP algorithm contains some more side details, such as problem preconditioning and solution polishing, for which we refer to \cite{OSQP:Stellato2020} for more information. See also Algorithm \ref{alg:osqp} in Section \ref{sec:benchmarks} for a compact OSQP summary.
\end{remark}

As long as the penalty parameter $\rho$ remains constant, the ADMM iterations are computationally inexpensive (consisting of a few matrix-vector multiplications and vector element-wise operations). Furthermore, for a static $\rho$, the convergence of the ADMM algorithm is guaranteed by showing that the following value function decreases along iterations, see \cite[Appendix A]{ADMMbook:Boyd2010} for details, i.e., 
\begin{equation}
    \bar{V}^k := (1/\rho)\|y^k-y^*\|_2^2 +\rho\|z^k-z^*\|_2^2,
\end{equation}
where $y^*$ and $z^*$ denote the optimal values for $y$ and $z$ with respect to Problem \eqref{eq:ADMM_prob_z}.
The convergence rate depends on the choice of relaxation parameter $\rho > 0$, which can be updated online to improve convergence. However, ADMM remains limited by a linear convergence rate. As a result, while iterations are inexpensive, it can take easily hundreds or more iterations to converge to a satisfying tolerance. In the next section, we show that by introducing a novel multi-dimensional weighting method, the convergence rate of the ADMM becomes dynamic and typically much faster.

\section{The \solver\ algorithm}\label{sec:superADMM}
The augmented Lagrangian in \eqref{eq:lagr_prox} illustrates that the penalty between $Ax$ and $z$ is weighted with a constant $\rho$. Notice that for each $z_i$ the optimal value is either $A_ix$ (unconstrained) or $l_i / u_i$ (clipped). This introduces a problem for the choice of $\rho$, i.e., for all unconstrained $z_i$, the preferred value $\rho$ is small. However, for clipped values of $z_i$, the optimal solution $x$ must be forced to the constraint, thus preferring a larger weight $\rho$. To circumvent this, we first weigh each constraint individually, by defining
\begin{equation}\label{eq:ADMM_rhotoR}
    \begin{aligned}
        R := \text{diag}\{\rho_1, \ldots, \rho_m\}
    \end{aligned}
\end{equation}
and we rewrite the augmented Lagrangian as
\begin{equation}\label{eq:lagranR}
    \begin{aligned}
        \mathcal{L}(x,x^k,z,y,R) :=& \frac{1}{2}x^TPx + x^Tq + g(z) + y^T(Ax-z) \\ &+ \frac{1}{2}\|Ax-z\|_R^2 + \frac{\sigma}{2}\|x-x^k\|_2^2,
    \end{aligned}
\end{equation}
where $\|\cdot\|_R$ denotes a weighted norm over $R\succ0$, i.e., $\|x\|_R =\sqrt{x^TRx}$. This yields the following ADMM steps
\begin{subequations}\label{eq:ADMM_algR}
    \begin{align}
    \label{eq:minxR} x^{k+1} &:= \text{arg}\min_x \mathcal{L}(x, x^k, z^k, y^k, R) \\
    \label{eq:minzR} z^{k+1} &:= \Pi_{\left[l, u\right]}(Ax^{k+1} + R^{-1}y^k),\\
    \label{eq:minyR} y^{k+1} &:= y^k +R(Ax^{k+1}-z^{k+1}).
    \end{align}
\end{subequations}
Next, instead of keeping $R$ static, we redefine $R^k = \text{diag}\{\rho_1^k, \ldots, \rho_m^k\}$ such that $R_{i,i}^k = \rho_i^k$ and update the diagonal entries after each ADMM iteration in \eqref{eq:ADMM_algR} as
\begin{equation}\label{eq:Rupdate}
    \begin{aligned}
        R_{i,i}^{k+1} &= \begin{cases}
            \alpha R_{i,i}^{k}, ~~~~ & \mathrm{if}~z^\kpp_i = l_i \lor z^\kpp_i = u_i \\
            (1/\alpha) R_{i,i}^k, ~~~~ & \mathrm{else}
        \end{cases},
    \end{aligned}
\end{equation}
for all $i=\{1, \ldots, m\}$, and some hyperparameter $\alpha > 1$. Using the Lagrangian in \eqref{eq:lagranR}, $x^{k+1}$ in \eqref{eq:minxR} is the solution of $\nabla_x\mathcal{L}(x,x^k,z^k,y^k,R) = 0$, i.e.,
\begin{equation}\label{eq:minx_deriv1}
    \begin{aligned}
        0 &= Px + q + A^Ty^k + A^TR(Ax-z^k) + \sigma(x-x^k).
    \end{aligned}
\end{equation}
We can rewrite \eqref{eq:minx_deriv1} more efficiently by defining $\nu = y^k + R(Ax - z^k)$, and substituting it in \eqref{eq:minx_deriv1},
\begin{equation}\label{eq:minx_deriv2}
    \begin{aligned}
        0 &= (P+\sigma I)x -\sigma x^k + q + A^T\nu \\
        0 &= R^{-1}(y^k-\nu) + Ax - z^k,
    \end{aligned}
\end{equation}
which is a symmetric system of equations in $x$ and $\nu$. Furthermore, if we obtain $x^\kpp$ and $\nu^\kpp$ as the solution of \eqref{eq:minx_deriv2}, we can obtain $Ax^\kpp$ also as $Ax^\kpp = z^k + R^{-1}(\nu^\kpp-y^k)$, which eliminates the need for an additional matrix-vector multiplication. Finally, using the Lagrangian in \eqref{eq:lagranR} and derivation \eqref{eq:minx_deriv2}, we can insert $R^{k}$, which is updated as in \eqref{eq:Rupdate}, which allows rewriting the full ADMM steps as follows:
\begin{subequations}\label{eq:ADMM_full}
    \begin{align}
        \label{eq:minx}\mat{ \sigma x^k- q \\ z^{k} - (R^k)^{-1}y^{k}}&=\mat{P+\sigma I & A^T \\ A & -(R^{k})^{-1}}\mat{x^{k+1}\\ \nu^{k+1}}  \\
        \label{eq:minzt}\tilde{z}^{k+1} &= z^k + (R^{k})^{-1}(\nu^{k+1} - y^k) \\
        \label{eq:minz} z^{k+1} &= \Pi_{\left[l, u\right]}(\tilde{z}^{k+1} + (R^{k})^{-1}y^k),\\
        \label{eq:miny}y^{k+1} &= y^k +R^{k}(\tilde{z}^{k+1}-z^{k+1}),
    \end{align}
\end{subequations}
where $\tilde{z}^\kpp = Ax^\kpp$.
\begin{remark}
    It is worth to point out that a dynamic weighting method for updating individual weights for each constraint via a diagonal matrix $R$ was originally employed in the QPALM solver \cite{hermans2019qpalm} (see also Algorithm \ref{alg:qpalm} in Section \ref{sec:benchmarks} for a compact summary of QPALM), however, with some notable differences. Therein, penalty parameters are also updated per constraint; however, as shown in line \ref{ln:qpalmupdate} in Algorithm \ref{alg:qpalm}, the update is made only if the residual of the corresponding constraint (i.e. $A_ix^k-z_i^k$) is not decreasing fast enough. Additionally, \solver\ also updates individual penalty parameters when constraints are not active and, furthermore, QPALM solves a different system of equations instead of \eqref{eq:minx}, see line \ref{ln:qpalmsolve} in Algorithm \ref{alg:qpalm} of Section \ref{sec:benchmarks}, which requires using a Newton method with line search.
\end{remark}

\subsection{Numerical Stability}
Exponentially increasing and decreasing the diagonal values in $R$ causes increasingly poor matrix conditioning, influencing the numerical precision of the $x^{k+1}$ updates. This prevents convergence of the algorithm (more specifically, the dual objective). To mitigate this issue, we limit $R^k$ as $1/b^k\leq R^k\leq b^k$ with some bound $b^k$, which is updated through
\begin{equation}\label{eq:ConditionNumberUpdate}
    \begin{aligned}
        b^{k+1} &= \begin{cases}
            \tau b^k, ~~ &\text{if } \varepsilon^{k+1} \geq r^{k+1}_\prim \\
            b^k, ~~ &\text{else} 
        \end{cases},
    \end{aligned}
\end{equation}
where
\begin{equation}\label{eq:ConditionNumber}
    \begin{aligned}
        \varepsilon^{k+1} &= \left\|\mat{\sigma x^{k} - q \\ z^{k} - (R^{k})^{-1}y^{k}} \right. \\ & ~~~~~~\left.-\mat{P+\sigma I & A^T \\ A & -(R^{k})^{-1}}\mat{x^{k+1}\\ \nu^{k+1}}\right\|_{\infty} \\ 
    \end{aligned}
\end{equation}
and $0<\tau<1$, which is set to $0.5$ as default.
Furthermore, if convergence to the desired accuracy is not reached before $b^k<1$, the algorithm cannot obtain a better solution (for the given parameters $\alpha$, $\tau$, and $\sigma$) and will thus end with an exit flag that indicates the solution is inaccurate. For this reason, \solver\ allows setting $\epsilon_\text{abs} = 0$, a case in which it runs until $b^k < 1$ or the iteration or time limits are exceeded. Initially, we set $b^0 = 10^8$ to limit the growth of $R^k$, yet due to the scale of the cost function, this can be changed in the solver settings.

The \solver\ algorithm based on the above derivations is formally summarized in Algorithm \ref{alg:solver}. The algorithm illustrates the standard ADMM iterates, between lines 6 and 9, and the novel dynamic weighting update equations in lines 11-12.
\begin{algorithm}[t]
    \caption{\solver}\label{alg:solver}
    \begin{algorithmic}[1]
    \Require $P,q,A,l,u$, $x_0$, $y_0$, $\rho_0, b^0, \tau$, $\alpha, \sigma$, $\epsilon_\mathrm{abs}$, $k_\mathrm{max}$
    \State Initialize $x^0 \gets x_0, y^0\gets y_0, z^0 \gets Ax_0$
    \State $R_{i,i}^0 \gets \rho_0, \forall i \in\{1,\ldots, m\}$
    \State $r^0_\prim = r^0_\dual = 2\epsilon_\mathrm{abs}$ \Comment{Ensures \textbf{while} loop starts}
    \State $k\gets 0$
    \While{ ($r_\prim^{k} > \epsilon_\text{abs}$ \textbf{or} $r_\dual^{k} > \epsilon_\text{abs}$) \textbf{and} $k<k_\mathrm{max}$}
        
        \State $x^{k+1}, \nu^\kpp \gets$ Solve \eqref{eq:minx} \label{ln:solve}
        \State $\tilde{z}^\kpp \gets z^k + (R^k)^{-1}(\nu^\kpp-y^k)$
        \State $z^{k+1} \gets \Pi_{\left[l, u\right]}(\tilde{z}^{k+1} + (R^k)^{-1}y^k)$
        \State $y^{k+1} \gets y^k +R^k(\tilde{z}^{k+1}-z^{k+1})$
        \State Compute $r_\prim^{k+1}$ and $r_\dual^{k+1}\gets$ \eqref{eq:convergenceTerms}
        \State Compute $\varepsilon^\kpp \gets$ \eqref{eq:ConditionNumber} and update $b^\kpp\gets $ \eqref{eq:ConditionNumberUpdate}
        \State $R_{i,i}^{k+1} = \begin{cases}
            \min(b^\kpp,\alpha R_{i,i}^{k}), ~\mathrm{if}~z^k_i = l_i \lor z^k_i = u_i \\
            \max(1/b^\kpp, (1/\alpha) R_{i,i}^k), ~\mathrm{else}
        \end{cases}$
        \State $k\gets k+1$
    \EndWhile
    \Ensure $x, y$
    \end{algorithmic}
\end{algorithm}

\begin{remark}
    The convergence of the developed \solver\ algorithm can be analyzed based on the results in \cite[Prop. 1-2]{ADMMproof:Bertsekas1976}. Therein, it has been shown that for an equality constrained problem and $\rho\rightarrow\infty$, the ADMM algorithm exhibits a superlinear convergence rate. Furthermore, in the proof, they consider $\rho = \alpha^k$, for any $\alpha>1$, as an example term. If problem \eqref{eq:ADMM_prob} would be equality constrained only, $R^k$ would behave in exactly the same way, which suggests that superlinear convergence can be formally proven for \solver. However, the complete convergence proof of the proposed method is beyond the scope of this paper and will be addressed in future work.
\end{remark}

\subsection{Infeasibility}
The infeasibility of the QP problem is checked following the method presented in \cite{OSQP:Stellato2020, ADMMinfeasibility:Banjac2019}. This states that for $\delta y^k = y^k-y^{k-1}$ if
\begin{equation}\label{eq:priminfeas}
    \begin{aligned}
        A^T\delta y^k &= 0, ~~ \text{and} ~~ u^T\delta y^k_+ + l^T\delta y^k_- < 0,
    \end{aligned}
\end{equation}
then the problem is primal infeasible. Similarly, given $\delta x^k = x^k - x^{k-1}$, if
\begin{equation}\label{eq:dualinfeas}
    \begin{aligned}
        q^T\delta x^k &< 0, ~~ P\delta x^k = 0, ~~ \text{and} \\
        A_i\delta x^k& \begin{cases} = 0 ~~ & l_i, u_i \in \mathbb{R} \\
            \geq 0 ~~ & u_i = +\infty, l_i\in\mathbb{R}\\
            \leq 0 ~~ & l_i = -\infty, u_i \in \mathbb{R},
        \end{cases}
    \end{aligned}
\end{equation}
the problem is dual infeasible. The infeasibility check is performed every 10 iterations to reduce unnecessary computation. Furthermore, \eqref{eq:priminfeas} and \eqref{eq:dualinfeas} also hold if $y^*$ and $x^*$ are unbounded, respectively \cite{ADMMinfeasibility:Banjac2019}.

\subsection{Parameter Selection}
The numerical stability and convergence speed of the proposed \solver\ algorithm is determined by the choice of parameters $\alpha$, $\sigma$ and $b^0$, which we discuss next. The parameter $\alpha > 1$ controls the rate of convergence. Setting the parameter too large can cause the solver to run into numerical issues at a too early state of convergence. We set it to $\alpha = 500$ for a decent balance. The parameter $\sigma$ should be kept small. However, higher values for $\sigma$ are desired if the problem is poorly scaled or ill-conditioned. We do not include adaptive procedures as done in \cite{hermans2019qpalm}, but keep it constant to $10^{-6}$, as done in \cite{OSQP:Stellato2020}.



\section{\Solver\ Implementation}\label{sec:implem}
The current implementation is written in \texttt{c} and supports dense and sparse problems. The dense formulation accelerates linear algebra computations with functions from CBLAS \cite{BLAS:Lawson1979} and LAPACK \cite{LAPACK}. The sparse solver relies on functions from \texttt{CSparse} \cite{CSparse:Davis2006} and uses a pivot-free sparse $LDL^T$ solver from \cite{LDLTsparse:Davis2005} to solve the linear system in Line \ref{ln:solve} of Algorithm \ref{alg:solver}. To reduce the number of fill-ins in $L$, we permute the matrix in \eqref{eq:minx} using the Approximate Minimal Degree (AMD) sorting method \cite{AMD:Davis1996}. Since the sparse structure of the matrix in \eqref{eq:minx} is static, we only perform the AMD and symbolic analysis (which determines the elimination tree of a matrix and counts the nonzeros per column) once. Because of this, we can allocate the necessary memory for iterative computations prior to the loop, which is suitable for memory control in embedded platforms. We developed wrappers for the function \texttt{\solver()} in MATLAB and Python. The current version of the solver is publicly available at: \texttt{https://github.com/Petrus1904/superADMM}.
\begin{remark}
    The \solver\ algorithm relies on double precision floating point calculations to achieve its fast convergence rate. Ensuring the solver also works effectively using single-precision float is of interest to implement it on embedded platforms, and remains part of future work.
\end{remark}
\begin{remark}
    Following \cite{OSQP:Stellato2020}, the updates for $x^{k+1}$ and $z^{k+1}$ can also be solved using the ``indirect method", which replaces \eqref{eq:minx} and \eqref{eq:minzt} to:
    \begin{equation}\label{eq:indirectMethod}
        \begin{aligned}
        x^{k+1} &= (P+ \sigma I + A^TR^{k}A)^{-1}\\ &  ~~~~~\times (A^T(R^{k}z^k - y^k)+\sigma x^k - q) \\
        \tilde{z}^{k+1} &= Ax^{k+1}
        \end{aligned}
    \end{equation}
    This method reduces the size of the inverse from $(n+m)^2$ to $n^2$. However, the indirect method also requires recomputing $(P+\sigma I+ A^TR^{k}A)$ every iteration (that is, the upper triangular part). For dense matrices, this is still faster than the direct method \eqref{eq:ADMM_full}.
\end{remark}


\begin{remark}
    The scale of the problem heavily influences the convergence of classical ADMM algorithms. Preconditioning the problem allows for standardizing default selections (i.e. $\rho^0$) and reduces the number of iterations. With \solver, we noticed that convergence is much less influenced by the scale of the original problem, making preconditioning often unnecessary, or can be fixed by reducing $\alpha$ or slightly increasing $\sigma$. However, this can become an option in future versions.
\end{remark}

\begin{remark}
The standard $LDL^T$ decomposition algorithm considers a diagonal $D$ matrix. As an alternative, we can implement the Bunch-Kaufman factorization method \cite{BunchKaufman:Bunch1977}, which changes some diagonal entries in $D$ to 2$\times$2 blocks. This improves the numerical conditioning of the decomposition and, as a result, the accuracy of $x^\kpp$. We have written a pivot-free Bunch-Kaufman factorization method on top of the standard pivot-free $LDL^T$ decomposition from \cite{LDLTsparse:Davis2005}. The results show that it can solve more ill-conditioned problems and it can solve the same problem in fewer iterations. However, this decomposition strategy results in a denser $L$ and is computationally more expensive. Therefore, while this approach solves problems in fewer iterations and is able to solve more complex problems without encountering numerical problems, to date it has not outperformed the standard algorithm $LDL^T$ in execution time. We aim to include this method as an option in future versions.
\end{remark}



\section{Comparison Benchmarks}\label{sec:benchmarks}
In this section, we benchmark the \solver\ solver against four modern off-the-shelf solvers, OSQP v0.6.2 \cite{OSQP:Stellato2020}, QPALM v1.2.3 \cite{hermans2019qpalm}, Mosek v9.3.20 \cite{mosek} and Gurobi v11.0 \cite{gurobi}. For the sake of comparison, we have outlined the OSQP and QPALM methods in Algorithms \ref{alg:osqp} and \ref{alg:qpalm}, respectively. For completeness, we define for Algorithm \ref{alg:qpalm}:
\begin{equation}
    \begin{aligned}
        \mathcal{J}(x) &:= \{i ~|~ (Ax + (R^k)^{-1}y^k)_i \notin [l_i, u_i]\} \\
        \nabla\phi_k(x) &:= Px + q + A^T(y^k+R^k(Ax-Z_k(x))) \\ &~~~~~~+\sigma^k(x-x^k),~~\text{where,} \\
        Z_k(x) &:= \Pi_{[l,u]}(Ax+(R^k)^{-1}y^k),
    \end{aligned}
\end{equation}
and $A_{\mathcal{J}(x)}$ denotes all rows in $A$ that belong to $\mathcal{J}(x)$. Please do mind that to align the OSQP and QPALM solvers with \solver, some notational changes can be observed between the algorithms presented here and the original works. Additionally, for conciseness, we left out the infeasibility checks (which are the same for all three solvers) and the relative termination condition, which is unused in this paper.

The benchmark consists of four different problem structures, most of which are taken from \cite{OSQP:Stellato2020}. For random distributions, we denote $\mathcal{N}(\text{mean}, \text{var})$ as a Gaussian distribution and $\mathcal{U}(\text{lb},\text{ub})$ as a uniform random signal.



\begin{algorithm}[t]
    \caption{OSQP \cite{OSQP:Stellato2020}}\label{alg:osqp}
    \begin{algorithmic}[1]
    \Require $P,q,A,l,u$, $x_0$, $y_0$, $\rho >0$, $\beta \in (0,2), \sigma >0$, $\epsilon_\mathrm{abs}$, $k_\mathrm{max}$
    \State Initialize $x^0 \gets x_0, y^0\gets y_0, z^0 \gets Ax_0$
    \State Precondition problem
    \State $r^0_\prim = r^0_\dual = 2\epsilon_\mathrm{abs}$ \Comment{Ensures \textbf{while} loop starts}
    \State $k\gets 0$
    \While{ ($r_\prim^{k} > \epsilon_\text{abs}$ \textbf{or} $r_\dual^{k} > \epsilon_\text{abs}$) \textbf{and} $k<k_\mathrm{max}$}
        
        \State Solve $\mat{\sigma x^k- q \\ z^{k} - \rho^{-1}y^{k}}=\mat{P+\sigma I & A^T \\ A & -\rho^{-1}I}\mat{\tilde{x}^\kpp\\ \nu^\kpp}$
        \State $\tilde{z}^\kpp \gets z^k + \rho^{-1}(\nu^\kpp-y^k)$
        \State $x^\kpp = \beta\tilde{x}^\kpp + (1-\beta)x^k$
        \State $z^{k+1} \gets \Pi_{\left[l, u\right]}(\beta\tilde{z}^{k+1} +(1-\beta)z^k + \rho^{-1}y^k)$
        \State $y^{k+1} \gets y^k +\rho(\beta\tilde{z}^{k+1}+(1-\beta)z^k-z^{k+1})$
        \State Compute $r_\prim^{k+1}$ and $r_\dual^{k+1}\gets$ \eqref{eq:convergenceTerms}
        \State $k\gets k+1$
    \EndWhile
    \Ensure $x, y$
    \end{algorithmic}
\end{algorithm}

\begin{algorithm}[t]
    \caption{QPALM \cite{hermans2019qpalm}}\label{alg:qpalm}
    \begin{algorithmic}[1]
    \Require $P,q,A,l,u$, $x_0$, $y_0$, $\alpha>1, b, \sigma^0$, $[\vartheta,\Delta_x, \Delta_\varepsilon]\in(0,1)^3$, $\varepsilon^0_\mathrm{in}, \epsilon_\mathrm{abs}$, $k_\mathrm{max}$
    \State Initialize $x^0 \gets x_0, y^0\gets y_0, z^0 \gets Ax_0$
    \State Precondition problem
    \State $r^0_\prim = r^0_\dual = 2\epsilon_\mathrm{abs}$ \Comment{Ensures \textbf{while} loop starts}
    \State $k\gets 0$
    \While{ ($r_\prim^{k} > \epsilon_\text{abs}$ \textbf{or} $r_\dual^{k} > \epsilon_\text{abs}$) \textbf{and} $k<k_\mathrm{max}$}
    \State $x\gets x^k$
        \While{$\sigma^k\|\nabla \phi_k(x)\|_2> \varepsilon^k_\mathrm{in}$}
            \State Solve \resizebox{0.70\columnwidth}{!}{$\mat{-\nabla\phi_k(x) \\ 0}=\mat{P+\sigma^k I & A_{\mathcal{J}(x)}^T \\ A_{\mathcal{J}(x)} & -(R^{k}_{\mathcal{J}(x)})^{-1}}\mat{d\\ \nu}$}\label{ln:qpalmsolve}
            \State Find optimal step size $t_\star$ using exact line search
            \State $x\gets x + t_\star d$
        \EndWhile
        \State $x^\kpp \gets x$, $z^{k+1} \gets \Pi_{\left[l, u\right]}(Ax^\kpp + (R^k)^{-1}y^k)$
        \State $y^{k+1} \gets y^k +R^k(Ax^\kpp-z^{k+1})$
        \State Compute $r_\prim^{k+1}$ and $r_\dual^{k+1}\gets$ \eqref{eq:convergenceTerms}
        \State \resizebox{0.90\columnwidth}{!}{$R_{i,i}^{k+1} = \begin{cases}R_{i,i}^k ~~ \mathrm{if} ~|(Ax^{k+1} - z^{k+1})_i| \leq \vartheta |(Ax^{k} - z^{k})_i|, \\ \min(\alpha R_{i,i}^k, b) ~~ \mathrm{else}\end{cases}$}\label{ln:qpalmupdate}
        \State $\sigma^\kpp \gets \Delta_x \sigma^k$
        \State $\varepsilon^\kpp_\mathrm{in} \gets \Delta_\varepsilon\varepsilon^k_\mathrm{in}$
        \State $k\gets k+1$
    \EndWhile
    \Ensure $x, y$
    \end{algorithmic}
\end{algorithm}

\subsection{Model Predictive Control}
Consider the following MPC problem:
\begin{equation}\label{eq:MPCprob}
    \begin{aligned}
        \min_{x,u} ~~&x_{N|k}^TQ_Tx_{N|k} + \sum_{i = 0}^{N-1} x_{i|k}^T&&\hspace{-0.35cm}Qx_{i|k} + u_{i|k}^TRu_{i|k} \\
        \text{s.t.}~~ &x_{0|k} = x_0 &&\\
        &x_{i+1|k} = \bar{A}x_{i|k} + \bar{B}u_{i|k},  & & \forall i \in \{0, \ldots, N-1\} \\
        &x_{i+1|k} \in \mathbb{X}, u_{i|k} \in \mathbb{U}, & & \forall i \in \{0, \ldots, N-1\}
    \end{aligned}
\end{equation}
where $x_{i|k}\in\mathbb{R}^{n_x}$ and $u_{i|k}\in\mathbb{R}^{n_u}$ denotes the predicted state and input, respectively, with $n_u = (1/2)n_x$. We generate the system dynamics as $A = I+\Delta$, with $\Delta_{ij} \sim\mathcal{N}(0,0.1)$, $B_{ij}\sim\mathcal{N}(0,1)$. $Q$ is diagonal with $Q_{ii}\sim\mathcal{U}(0,10)$ and $R = 0.1I_{n_u}$. The terminal cost $Q_T$ is the solution of the discrete algebraic Riccati equation such that a stabilizing MPC control law is obtained. The state and input are constrained with sets $\mathbb{X} := \{ x ~|~ \ubar{x}\leq x \leq \bar{x}\}$, with $\ubar{x} = -\bar{x}$ and $\bar{x}_i \sim\mathcal{U}(1,5)$. The set $\mathbb{U}$ is designed similarly, also with $\bar{u}_i \sim\mathcal{U}(1,5)$. The control goal is to steer the initial state $x_0\sim\mathcal{U}(0.5\ubar{x}, 0.5\bar{x})$ to zero. Two MPC benchmarks have been included. First, we fix the prediction horizon $N=10$ and vary the state order $n_x$, secondly, we fix the state order at $n_x=10$ and the prediction horizon is varied.

\subsection{Data fitting: Lasso and Huber loss}
Finally, we compare different robust data fitting methods, each designed to handle data outliers better than least--squares. \textit{Lasso} is a linear regression technique that weights the optimization variable with an $l_1$ norm, i.e.,
\begin{equation}
    \begin{aligned}
        \min_{x} \|Ax-b\|^2_2 + \lambda \|x\|_1.
    \end{aligned}
\end{equation}
where $x\in\mathbb{R}^{n}$ and $A\in\mathbb{R}^m$, with $m=100n$. This can be converted to a QP as \cite{OSQP:Stellato2020}:
\begin{equation}
    \begin{aligned}
        \min_{y,x,t} & ~~y^Ty + \lambda\mathbf{1}^Tt \\
        \text{s.t.} & ~~ y=Ax-b \\
        &~~ -t\leq x\leq t,
    \end{aligned}
\end{equation}
where we generate $A_{ij}\sim\mathcal{N}(0,1)$ with 15\% nonzero elements. The data in $b$ is generated as $b=Av+w$, where
\begin{equation}
    v \sim \begin{cases}
        0 & \text{with probability } p=0.5 \\
        \mathcal{N}(0,1/n) & \text{otherwise,}
    \end{cases}
\end{equation}
and $w_i \sim\mathcal{N}(0,1)$. Finally, as done in \cite{OSQP:Stellato2020}, $\lambda = \frac{1}{5}\|A^Tb\|_{\infty}$.

\textit{Huber loss} is a robust regularization technique for data with outliers. This results in the following QP \cite{OSQP:Stellato2020}:
\begin{equation}
    \begin{aligned}
        \min_{x,u,r,s} ~~& u^Tu + 2M\mathbf{1}^T(r+s) \\ 
        \text{s.t.} ~~& Ax-b-u = r-s \\
         & r\geq 0, ~~ s\geq 0.
    \end{aligned}
\end{equation}
We generate $A\in\mathbb{R}^{n\times m}$, with $m=100n$ and $A_{ij}\sim\mathcal{N}(0,1)$ with 15\% nonzero elements. The data in $b$ is generated as $b=Av+w$, where $v_i \sim\mathcal{N}(0,1/n)$ and \cite{OSQP:Stellato2020}
\begin{equation}
    w_i \sim \begin{cases}
        \mathcal{N}(0,1/4) & \text{with probability } p=0.95 \\
        \mathcal{U}(0,10) & \text{otherwise.}
    \end{cases}
\end{equation}

\subsection{Results}\label{ssec:results}
For all solvers, we set the solution tolerance $\epsilon_\text{abs} = 10^{-8}$ and the relative tolerance in OSQP and QPALM $\epsilon_\text{rel} = 0$ as this is currently not supported in \solver. All benchmarks are tested over a range of problem sizes $n$, each tested five times for consistency. For all timing results, we considered the internal runtime as reported by each solver in the results. However, Mosek does not report any runtime that is less than 0.015s (it reports 0.0s). For this reason, we measure the total time in MATLAB and use that if Mosek fails to report a nonzero runtime. Something similar happens with Gurobi for an execution time less than 0.005s, however, it still shows the work time. Please note that, therefore, the results of Mosek and Gurobi are offset for smaller problems.

\begin{figure}[h]
	\centering
	\includegraphics[width=1\columnwidth, trim={0.5cm, 0cm, 2.0cm, 0cm}, clip]{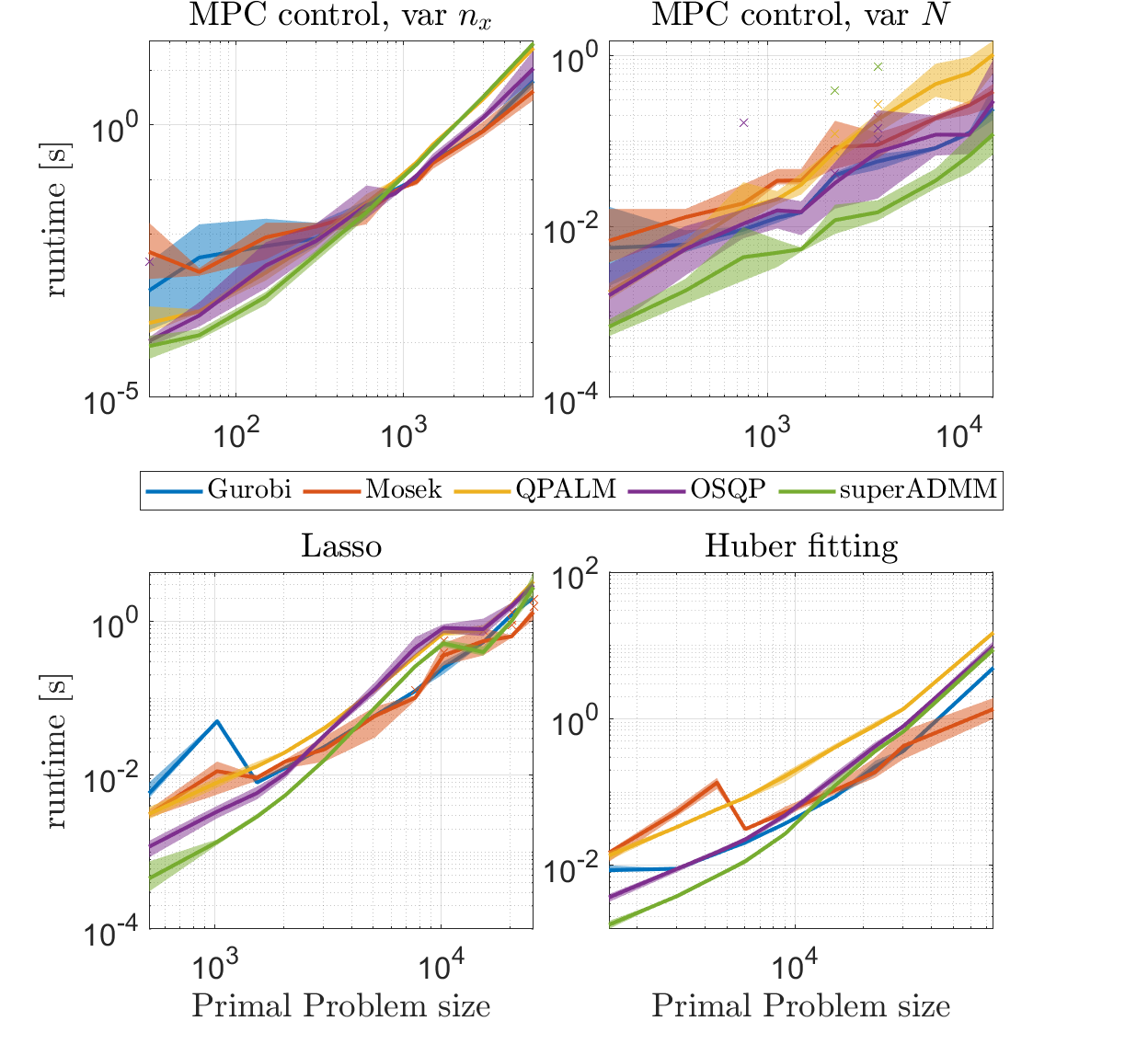}
	\caption{Timing results corresponding to Gurobi, Mosek, OSQP, QPALM and \solver\ for different benchmark problems.}
	\label{fig:benchmarks}
\end{figure}

All experiments were carried out using MATLAB on a laptop with an Intel i7-9750H @2.60GHz CPU and 16GB of RAM. The results are displayed in Figure \ref{fig:benchmarks}. The solid line shows the mean, while the semi-transparent patch displays the minimum and maximum run times. Every cross ($\times$) indicates a problem the solver failed to solve (infeasible or ran out of iterations) or solved inaccurate.

For small to medium scale problems, \solver\ outperforms all the other tested solvers for all 4 benchmarks. Additionally, \solver\ shows very consistent results, with a very small variation for all problem sizes, which is a nice feature for real-time implementation.
As the problems become large scale, Mosek is the most efficient solver. Interestingly, while \solver\ is outperformed by all other solvers for large-scale MPC problems with an increasing state-dimension, \solver\ was the best performing solver for MPC problems for increasing prediction horizons. For Lasso and Huber benchmarks, \solver\ is only outperformed by Gurobi and Mosek for large scale problems, which is expected, given they use an interior-point method.

\subsection{Maros-Mészáros Benchmarks}\label{ssec:MMbench}
To assess the solver's performance on challenging or ill-conditioned problems, we additionally let it solve the Maros-Mészáros benchmarks \cite{MarosMeszarosTests:1999, qpbenchmark2024}. This set includes 138 complex problems on diverse scales. To execute the tests as fair as possible, we set the maximum iterations of all solvers to infinity, the time limit to $t_\mathrm{lim} = 120$~seconds and $\epsilon_\mathrm{rel} = 0$ for OSQP and QPALM. Table \ref{tbl:MM} presents the number of problems successfully solved by each solver at a low tolerance level ($\epsilon_\mathrm{abs} = 10^{-3}$) versus a high tolerance level ($\epsilon_\mathrm{abs} = 10^{-9}$). The results show that QPALM solves 8 problems more than \solver\ in the case of $10^{-3}$ precision, while \solver\ solves 9 problems more than QPALM in the case of $10^{-9}$ precision, which shows the benefit of the novel dynamic weighting update method employed by \solver\ in handling ill-conditioned problems. Both QPALM and \solver\ solved more problems than OSQP for low- and high-precision settings. Table~\ref{tbl:MM} also shows the number of problems that have not been solved before the time limit was exceeded. For high tolerance benchmarks, 11 of these problems have not been solved by any solver; for low tolerance benchmarks, this was just two (\texttt{CONT-300} and \texttt{CVXQP3\_L}).
\begin{table}[h!]
\centering
 \begin{tabular}{| l l | c c c|} 
 \hline
 $\epsilon_\mathrm{abs}$ & & OSQP & QPALM & \solver  \\ [0.5ex] 
 \hline 
\multirow{2}{*}{$10^{-3}$} & solved & 114 & 131 & 123 \\
& $t_\mathrm{run} > t_\mathrm{lim}$ & 22 & 5 & 5 \\
\hline
\multirow{2}{*}{$10^{-9}$} & solved & 77 & 87 & 95 \\  
& $t_\mathrm{run} > t_\mathrm{lim}$ & 60 & 47 & 16 \\[0.5ex] 
 \hline
 \end{tabular}
 \caption{Number of problems solved of the Maros-Meszaros Benchmarks}
 \label{tbl:MM}
\end{table}

\subsection{Discussion}
The \solver\ solver shows competitive benchmark results, comparing well in execution speed and handling ill-conditioned problems. The results in subsection \ref{ssec:results} indicate that \solver\ is most effective at handling small- and medium-sized problems. However, it is less efficient than the interior-point methods used in Gurobi and Mosek for large-scale problems. In fact, for the MPC benchmark with varying $n_x$, it was outperformed by all solvers for the largest-scale problem (which is at $n_x = 400$, so 6000 optimization variables). Our detailed timing analysis shows that over 98\% of the 30 seconds spent solving these problems is spent decomposing and solving \eqref{eq:minx}. This indicates that solving the system of equations in \eqref{eq:minx} using the pivot-free $LDL^T$ method is specifically time-consuming for sparse matrices with over 600 nonzero values per column. However, for the MPC benchmark with increasing $N$, \solver\ outperformed the state-of-the-art regardless of size. This shows that with \solver, real-time MPC control of general linear systems can be achieved faster, or at the same time, but with a larger prediction horizon. 

The results of OSQP and QPALM in Table~\ref{tbl:MM} show that most of the unsolved problems are due to the time limit. This is commonly caused by the convergence rate becoming increasingly small as a result of poor numerical scaling. The forced convergence of \solver\ instead prevents this from happening, and if it does, we end up with a ``solved inaccurate" exitflag. From the unsolved problems in \solver, this was the major failure reason for termination.  

\subsection{Assessment of superlinear convergence}
\begin{figure}[h]
	\centering
	\includegraphics[width=0.90\columnwidth, trim={0.5cm, 0cm, 0.5cm, 0cm}, clip]{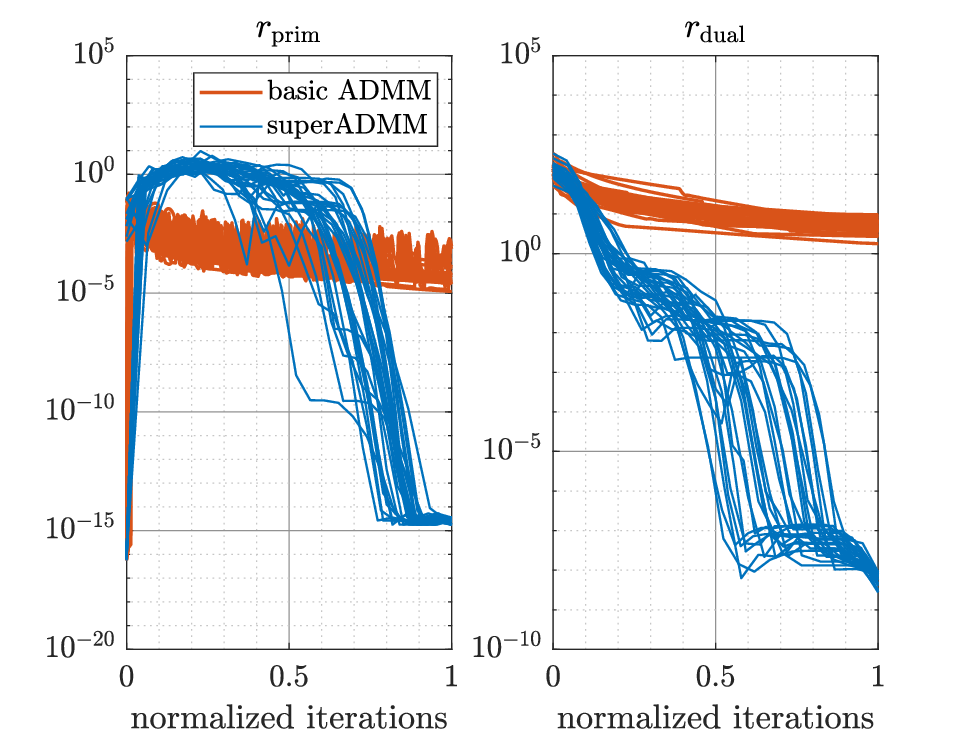}
	\caption{$r_\prim$ and $r_\dual$ over the normalized iteration (i.e., $\frac{k}{k_\mathrm{term}}$) for basic ADMM and \solver, with $\alpha = 10$ for visibility.}
	\label{fig:superlinear}
\end{figure}

In this paper, we have mentioned the superlinear convergence of \solver. In this subsection, we shall illustrate the effect of superlinear convergence due to the novel dynamic weighting method on the penalty matrix $R$. For this, we tested 30 random problems of medium size ($200$ primal variables, $300$ constraints). Figure \ref{fig:superlinear} displays the primal and dual residuals on a normalized iteration scale (that is, the current iteration divided by the termination iterate $k_\mathrm{term}$). For the basic ADMM, we also used \solver\ but set $\alpha = 1$, which disables dynamic weighting, thus obtaining standard ADMM. Here, we immediately see that the basic ADMM did not terminate within the maximum number of iterations, now set at 300. In contrast, \solver\ solved all problems within an average of 22 iterations (maximum of 29). 

The superlinear convergence of \solver\ is visible in Figure \ref{fig:superlinear} by the iteratively increasing convergence rate a few iterations before terminating. As the plots are logarithmic in the $y$-scale, it must be noted that this ``waterfall'' effect causes the primal residuals to reduce from $10^{-5}$ to $10^{-10}$ in a few iterations.

\section{conclusions}\label{sec:conclusions}
In this paper, we proposed a novel method to solve quadratic programs: \solver.  Unlike regular ADMM solvers, \solver\ employs a varying weighting method that weights each constraint individually and increases these weights at an exponential rate. Thanks to this, \solver\ can solve problems with high solution accuracy in just a few iterations. The method has been written in \texttt{c} with efficient linear algebra packages to provide a fast execution and can be invoked from popular languages like MATLAB and Python.

\bibliographystyle{IEEEtran}
\bibliography{references}
\end{document}